\documentclass[11pt,letterpaper]{article}
\usepackage{graphicx}

\newtheorem{THM}{Theorem}

\newtheorem{theorem}[THM]{Theorem}
\newtheorem{prop}[THM]{Proposition}

\newtheorem{lemma}[THM]{Lemma}

\newenvironment{Proof}{\textsc{Proof.}}{\QED}
\newcommand{\QED}{\hspace{8mm}\mbox{\textsc{qed}}\smallskip}

\newtheorem{claim}{Claim}

\newcommand{\curlyP}{\mathcal{P}}
\newcommand{\chiP}{\chi_\curlyP}
\newcommand{\minP}{m_\curlyP}

\newcommand{\chiD}{\chi_{\rm{d}}}
\newcommand{\tdc}{\chi_{\td}}
\newcommand{\td}{\mathrm{td}}
\renewcommand{\int}{\mathop{\mathit{int}}}

\newcommand{\un}{\mathcal{U}}
\newcommand{\conndom}{\mathop{\gamma_c}}

\newcommand{\property}[1]{\mbox{\textsf{#1}}}
\newcommand{\Dom}{\property{Dom}}
\newcommand{\TDom}{\property{TDom}}
\newcommand{\IF}{\property{IF}}
\newcommand{\Pone}{\property{Edge}}
\newcommand{\CDom}{\property{CDom}}
\newcommand{\Ptwo}{\property{Connected}}

\begin{document}

\begin{center}  {\large 
\textbf{Compelling Colorings: A generalization of the dominator chromatic number} \\[2mm] }
Anna Bachstein$^1$, Wayne Goddard$^{1,2}$, Michael A. Henning$^2$, John Xue$^{1}$ \\[2mm]
\textsuperscript{1}{School of Mathematical and Statistical Sciences, Clemson University} \\
\textsuperscript{2}{Dept of Mathematics and Applied Mathematics, University of Johannesburg} \\
\end{center}

\begin{abstract}
We define a $\curlyP$-compelling coloring as a proper coloring of the vertices of a graph such that every subset consisting of one vertex of each color has property $\curlyP$.
The  $\curlyP$-compelling chromatic number is the minimum number of colors in such a coloring. We show that this notion generalizes the dominator
and total dominator chromatic numbers, and provide some general bounds and algorithmic results. We also investigate
the specific cases where  $\curlyP$ is that the subset contains at least one edge or that the subset is connected.
\end{abstract}

\section{Introduction}

This paper generalizes the notion of the dominator chromatic number and the total dominator chromatic number of a graph. 

Recall that a set $S$ of vertices is a \emph{dominating set} of a graph $G$ if every vertex outside $S$ has a neighbor in $S$.
Further, a vertex/subset \emph{dominates} a set $T$ if it dominates every vertex in $T$. A \emph{dominator coloring} of a graph $G$ is a proper coloring of $G$ with the additional property that every vertex in the vertex set $V(G)$ dominates all vertices in at least one color class; that is: each vertex of the graph belongs to a singleton color class or is adjacent to every vertex of some (other) color class. The \emph{dominator chromatic number} $\chiD(G)$ of $G$ is the minimum number of color classes in a dominator coloring of $G$. This concept was introduced in \cite{CHH79} and \cite{HHMRH} and studied in \cite{BC12,Gera07a,Gera07b} inter alia.

The total version of dominator coloring is defined analogously. 
A set $S$ of vertices is a \emph{total dominating set} of a graph $G$ with no isolated vertices if every vertex in $V(G)$ has a neighbor in $S$.
A \emph{total dominator coloring} of $G$ is a proper coloring of $G$ in which each vertex of the graph is adjacent to every vertex of some other color class (different from its own color class). The \emph{total dominator chromatic number} of $G$, denoted by $\tdc(G)$ or $\chi_{d}^t(G)$, is the minimum number of color classes in such a coloring. This concept was introduced in~\cite{HHHMR} and 
studied in~\cite{Henning15,Kazemi15,Vijayalekshmi12} inter alia.

In this paper we consider another way of looking at these parameters.
Consider a graph $G$ where the vertices are properly colored.
Define a \emph{rainbow committee} (RC) as a set consisting of one vertex of each color.
Let $\curlyP$ denote some property of subsets of vertices: for example, this property might be that the set is dominating, forms a clique, or is connected.
We say that:
\begin{quote}
 a proper coloring of $G$ \emph{compels $\curlyP$} if every RC of $G$ has property $\curlyP$. 
 \end{quote}
Such a coloring might not exist, but if it does, we 
define the \emph{$\curlyP$-compelling chromatic number} of the graph $G$ as the minimum number of colors in a proper coloring that compels property $\curlyP$. We denote this by $\chiP(G)$.

As an example, consider the $5$-cycle. Up to symmetry there is a unique  coloring using $3$ colors. Here an RC consists either of three consecutive vertices or of two consecutive vertices and one isolated vertex. Thus every RC is dominating, but not all are connected. There is also an essentially unique coloring using $4$ colors; the resultant RC is connected. Thus the $\curlyP$-compelling chromatic number for $C_5$ is $3$ where $\curlyP$ represents dominating, but $4$ where  $\curlyP$ represents being connected.

In this paper we investigate the  $\curlyP$-compelling chromatic number for various properties $\curlyP$.
We start by showing that the (total) dominator chromatic number is equivalent to specific choices of $\curlyP$.
We then provide some general bounds on compelling numbers.
Thereafter we investigate two specific properties of the RC: namely, having at least one edge and being connected. 
We conclude with some comments about algorithmic questions.

\section{Some Equivalences}

We define the following properties of the rainbow committee:

\begin{quote}
Property \Dom: every vertex in the graph is dominated by the RC.\\[2mm]
Property \TDom: every vertex in the graph has a neighbor in the RC.
\end{quote}

\begin{theorem}
For all graphs $G$, \\
(a) 
the dominator chromatic number equals the \Dom-compelling chromatic number. \\
(b) 
the total dominator chromatic number equals the \TDom-compelling chromatic number. 
\end{theorem}
\begin{Proof}
We prove the result for the domination case; the proof for the total domination case is similar.
It suffices to show that a proper coloring is a dominator coloring if and only if it  is \Dom-compelling.

Consider a dominator coloring. Let $J$ be an RC and
consider any vertex~$u$. Since the vertex $u$ dominates some color class, say red,
and $J$ contains a red vertex, it follows that $J$ contains at least one vertex from the  closed neighborhood of $u$. Thus the coloring is \Dom-compelling.
On the other hand, consider a proper coloring that is not a dominator coloring.
Then there is some vertex $v$ that dominates no color class. Thus one can 
choose for each color a vertex not in $v$'s closed neighborhood to yield an RC that is not dominating. Thus the coloring is not \Dom-forcing.
\end{Proof}

Another natural property is that the subgraph induced by the rainbow committee has no isolates:

\begin{quote}
Property \IF: every vertex in the RC has a neighbor in the RC.
\end{quote}

But it turns out that the \IF-compelling chromatic number is also equivalent to the total dominator chromatic number.

\begin{theorem}
A proper coloring of a graph $G$ is \IF-compelling if and only if it is \TDom-compelling.
\end{theorem}
\begin{Proof}
It is immediate that any coloring that compels an RC to be totally dominating
also compels it to be isolate-free.

Suppose there is a proper coloring of $G$ that is \IF-compelling but not \TDom-compelling.
Then there is an RC $J$ where some vertex $v$ of the graph has no neighbor in $J$.
If $v$ is not in the set $J$, then one can change $J$ by choosing~$v$ instead of the
other vertex of its color. That is, we may assume that $v$ is in the~RC $J$. But then $v$ is isolated in $J$, which is
a contradiction, since we assumed that Property~\IF{} is compelled.
\end{Proof}


\section{Some General Bounds}

A relevant parameter is the minimum cardinality of a subset of $V(G)$ with property $\curlyP$. We denote this by $\minP(G)$.
Recall that $\chi(G)$ is the chromatic number of the graph $G$.
The following result generalizes the results for dominator colorings~\cite{Gera07a,Gera07b} and total dominator colorings~\cite{Kazemi15,Vijayalekshmi12}.
We say a property $\curlyP$ of subsets of $V(G)$  is \emph{upwards-closed} if every superset of a set with the property also has the property.
For example, ``is a dominating set of $G$'' or ``contains an edge'' is upwards-closed; ``induces a connected subgraph'' is not.

\begin{prop}  \label{p:general}
Let $G$ be a graph and consider a property $\curlyP$ of subsets of $V(G)$.\\
(a) $ \chiP (G) \ge \max \{ \, \minP (G) , \chi  (G)  \, \} $ \\
(b) If $\curlyP$ is upwards-closed then $\chiP (G) \le  \minP(G) + \chi(G)$.
\end{prop}
\begin{Proof}
(a) The first lower bound allows the RC to have property $\curlyP$.
The second lower bound allows a proper coloring.

(b) The upper bound is provided as follows: take a minimum subset $U$ with property $\curlyP$ and give each vertex of $U$ a unique
color; then color the remainder of the graph with a minimum proper coloring using new colors.\end{Proof}

These bounds can be achieved. For example,
for $G$ a complete bipartite graph it holds that $\chiP( G ) = \minP(G)$ provided 
$\minP(G)$ is achieved by a subset containing vertices from both partite sets. (For coloring: take a 
minimum subset with property $\curlyP$ and give each vertex a different color and then extend in any way to a proper coloring.)
On the other hand, the complete graph has only one proper coloring and so $\chiP(K_n) = \chi(K_n)$ if
$V(G)$ has the property.


We conclude this section with some comments on the disjoint union.
We say a property $\curlyP$ \emph{distributes over disjoint union} if for all 
subsets $S_1$ of $V(G_1)$ and subsets $S_2$ of $V(G_2)$, the union 
$S_1 \cup S_2$ has property~$\curlyP$ in the disjoint union of $G_1$ and $G_2$
if and only if $S_1$ has property $\curlyP$ in $G_1$ and $S_2$ has property $\curlyP$ in $G_2$.
For example, ``induces a bipartite subgraph'' distributes over disjoint union but ``contains an edge'' or ``induces a connected subgraph'' does not.

The following proposition is akin to the results for dominator coloring~\cite{GHR06} and total dominator coloring~\cite{Vijayalekshmi12}:

\begin{prop}
\label{p:disjointUnion}
Consider a property $\curlyP$ that distributes over disjoint union.
If $G$ is a disconnected graph with components $G_1, \ldots, G_k$, then 
\[
   \max_i    \left(  \chiP( G_i ) + \sum_{j\neq i}  \minP ( G_j )  \right) 
   \le \chiP( G ) 
   \le  \sum_{i}  \chiP( G_i ) \,.
\]
\end{prop}
\begin{Proof}
Consider a valid coloring of the graph $G$ and some RC~$J$.
For each component $G_j$ of $G$ let $u_j$ denote the number of colors that appear
only on~$G_j$. One must be compelled to choose a subset from 
that component $G_j$ with property $\curlyP$; so 
 $u_j \ge \minP(G_j)$. 
Further, since the restriction of the coloring to each component is a valid coloring thereof, for each $j$
there must be at least $\chiP(G_j)$ colors that appear on $G_j$. Thus 
\[
   \chiP( G ) \ge  \sum_j u_j  \, + \,  \max_i \, (\chiP(G_i)-u_i)  .
\]
The lower bounds follows.

The upper bound follows from using a disjoint palette for each component.\end{Proof}

There is equality throughout the chain in Proposition~\ref{p:disjointUnion} if every component $G_i$ has $\minP (G_i) = \chiP (G_i)$.

\section{Edge-compelling Chromatic Number}

In this section we consider the property that the RC is not an independent set.
Define

\begin{quote} 
Property \Pone:
\textit{at least one pair of vertices of the RC are joined by an edge}
\end{quote}

By Proposition~\ref{p:general} we get the following bound.
(The upper bound is equivalent to selecting one edge and giving its ends unique colors and giving the remainder of the graph any proper coloring.)

\begin{prop}  \label{p:generalP1}
For a graph $G$, the \Pone-compelling chromatic number is at least $\chi(G)$ and at most $\chi (G)+ 2$.
\end{prop}

\subsection{Paths and cycles}

\begin{theorem}
\label{t:p1path}
The \Pone-compelling chromatic number for the path $P_n$ is: \\
$2$ for $n=2,3$; \
$3$ for $n=4,5,6$; and \
$4$ for $n\ge 7$.
\end{theorem}
\begin{Proof}
It is immediate that the bipartite coloring is compelling exactly when $n=2,3$. A compelling coloring with three colors is not hard to find for $n=4,5,6$: here is 
one possibility for $n=6$:

\begin{center}
\includegraphics{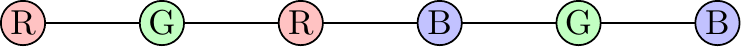}
\end{center}

By Proposition~\ref{p:generalP1}, the \Pone-compelling chromatic number of a path is at most $4$. We show next that 
the $7$-path requires $4$ colors.
Consider a proper coloring of the path $v_1\ldots v_7$ with three colors. Say vertex $v_4$ has color red.
Some other color must appear in $\{v_1,v_2\}$, say blue.
If one of $\{v_6,v_7\}$ has the third color then we are done, since one can choose an RC that is independent.
So we may assume the pair $\{v_6,v_7\}$ has colors red and blue. Similarly the pair $\{v_1,v_2\}$ has colors red and blue.
Now we may assume the third color appears at vertex $v_3$. Choose vertices $v_1$, $v_3$ and one of $\{v_6,v_7\}$ and 
one gets an RC that is an independent set, a contradiction.

Finally, consider the $n$-path for $n\ge 8$ and 
suppose one can compel with three colors.
Consider the first $7$ vertices: call this subpath $Q$. If $Q$ receives only $2$ colors,
then take for the RC of $P_n$ the first and fourth vertex of $P_n$ and any vertex of the third color.
If $Q$ has $3$ colors, then by the above discussion it is possible to find an RC from
this copy that is independent. In each case one obtains a contradiction.
\end{Proof}

\begin{theorem}
The \Pone-compelling chromatic number for the cycle $C_n$ is: \\
$2$ for $n=4$; \
$3$ for $n=3,5,6,7$; and \
$4$ for $n\ge 8$.
\end{theorem}
\begin{Proof}
It is easy to see that the bipartite coloring of only the $4$-cycle is compelling. So all other cycles need at least three colors.
Any coloring of the $3$- and $5$-cycle with three colors is compelling. A compelling coloring for the $6$-cycle is obtained
by taking the compelling coloring for the $6$-path shown above and joining the ends; 
and a compelling coloring for the $7$-cycle  is given in the following figure.
It follows that the compelling number of the $n$-cycle for $n=3,5,6,7$ is $3$.

\begin{center}
\includegraphics{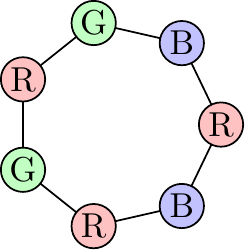} 
\end{center}

For $n\ge 8$, a compelling coloring using four colors is as follows: color two adjacent vertices
with unique colors and then color the remaining path with the other two colors.

It remains to show that three colors is not compelling for $n\ge 8$.
Take the colored $n$-cycle and discard some vertex that has a duplicate color.
What remains is a copy of the $(n-1)$-path with $3$ colors, where by Theorem~\ref{t:p1path} 
one can find an~RC that is independent.
\end{Proof}

\subsection{Graphs with small compelling number}

The \Pone-compelling chromatic number equals $2$ if and only if the graph is complete bipartite.
For, to equal $2$, one needs that the graph is bipartite. The two vertices in every RC must be adjacent: so the graph must be
complete bipartite.

By Proposition~\ref{p:generalP1} the compelling number is 
at most $4$ for trees. So when is it~$3$?
Recall that a graph has radius at most $2$ if there is some vertex $v$ such
that every vertex is within distance $2$ of $v$. We define a \emph{double-broom} 
as the caterpillar obtained by taking two stars $S$ and $S'$ and joining one leaf of
$S$ to one leaf of $S'$. See below for an example.

\begin{center}
\includegraphics{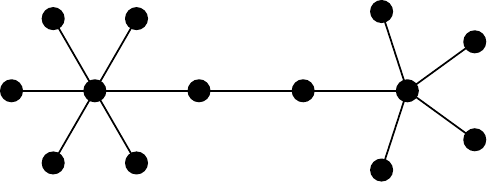}
\end{center}

\begin{theorem}
A tree $T$ has \Pone-compelling chromatic number $3$ if and only if $T$
has radius $2$ or is a double broom.
\end{theorem}
\begin{Proof}
We show first the sufficiency.
Assume the tree $T$ has radius~$2$, and let $v$ be a central vertex. Give $v$ 
the first color, its neighbors the second color, and the remaining vertices
the third color. This is a proper coloring.
Further, every RC contains both $v$ and a neighbor of $v$. That is, the coloring is \Pone-compelling.

For a double-broom, use the 
first color on only the center of one star $S$ and color the remainder of the tree 
properly using the remaining two colors. Every RC contains
the center of $S$; if it contains a second vertex of $S$ we are done. But 
otherwise it contains two vertices from $S'$ of different color, which must therefore be adjacent. That is, the coloring is \Pone-compelling.

So it remains to show that any other tree $T$ fails. 
A tree of radius $1$ is $K_2$ and that has parameter $2$.
Assume the radius of $T$ is larger than $2$. Then the diameter is at least $5$. Let $Q$ be a longest path in $T$.
Assume $Q$ has at least seven vertices. If $Q$ contains a $3$-colored copy of $P_7$ then 
by the above result on paths, $Q$ contains an RC without an edge. On the other hand if $Q$ is $2$-colored,
then it is easy to form an independent RC by starting with any vertex of the third color and adding two different colored vertices of $Q$.

So assume $Q$ has six vertices, say $v_1v_2 \ldots v_6$. By the maximality of $Q$, every neighbor of 
$v_2$ except $v_3$ is a leaf, and similarly every neighbor of $v_5$ except $v_4$ is a leaf.
We are assuming the tree is not a double broom. 
So there is some vertex~$w$ adjacent to one of the central vertices, say $v_3$. Say $w$ has color red. Then 
at least one of $X=\{v_1,v_2\}$ has a color other that red, say blue. 
If the third color appears on any of $Y=\{v_4,v_5,v_6\}$, then one obtains an RC that is independent by choosing it, $w$, and the blue vertex of $X$.
So all vertices on $Y$ are red or blue. By symmetric reasoning the third color cannot appear on~$X$, else one obtains an RC by taking the relevant
vertex, $w$ and a blue vertex of $Y$. Let $s$ be a vertex with the third color. It cannot be adjacent to both $v_2$ and $v_5$. So assume it is
not adjacent to $v_2$. Then we can take $s$, $v_6$,  and whichever of $\{v_1,v_2\}$ has color different to $v_6$ to obtain an independent RC. A contradiction.
\end{Proof}

This poses the question of what, in general, do graphs look like that have compelling number $3$.
For example the positive case of the above theorem generalizes to any triangle-free graph with a pair of nonadjacent vertices that dominate, such as
the $3$-cube. The negative result also generalizes:

\begin{theorem}
A graph $G$ with \Pone-compelling chromatic number $3$ has diameter at most $5$.
\end{theorem}
\begin{Proof}
Assume the diameter of $G$ is at least $6$. That is, there exist vertices
$v_1$ and $v_7$ joined by a shortest path $Q$ with internal vertices $v_2, \ldots, v_6$.
Note that $Q$ is induced.
If the  path $Q$ contains all three colors,
then we know by Theorem~\ref{t:p1path} that it contains an independent RC triple.
So we may assume that the vertices on $Q$ alternate colors, say red and blue. 

Consider a vertex $g$ of the third color. Since $Q$ is a shortest path, 
the neighborhood of $g$ is contained within some three consecutive vertices on $Q$. 
It is easily checked that one can find a red vertex $r$ and a blue vertex $b$ on $Q$ that are not
adjacent to each other nor to $g$. So $\{r,b,g\}$ is an independent RC, a contradiction.
\end{Proof}

Note that a graph with \Pone-compelling chromatic number $3$ can contain arbitrarily long induced paths: the fan that is the join of a dominator and a path is such an example.
We next show that 
a graph with \Pone-compelling chromatic number $4$ and large diameter is ``almost bipartite'' in some sense.

\begin{theorem}
There is a positive integer $D$ such that 
every graph $G$ with \Pone-compelling chromatic number $4$ and diameter at least $D$ contains
a complete bipartite subgraph $H$ such that $G-H$ is bipartite.
\end{theorem}
\begin{Proof}
We make no attempt to find the least diameter that works.
Assume the diameter is at least 73 and let $Q$ be a diametrical path. 
Then there are 19 pairs of consecutive vertices on $Q$ with 
at least two vertices between each pair (the 1st and 2nd vertex, the 5th and 6th vertex, etc).
Consider a compelling coloring. 
By the pigeon-hole principle, there must be four pairs that use the same two colors, say
orange and purple. Every vertex can be adjacent to at most one pair.
Thus, given any two vertices of the other colors, say red vertex $r$ and blue vertex $b$, 
one can find two pairs no vertex of which is adjacent to $r$ or $b$; and by taking an
orange vertex $o$ on one pair and a purple vertex $p$ on the other pair, one has a RC
where the only possible adjacency is between $r$ and $b$.
That is, one can find an independent RC unless $r$ and $b$ are adjacent for all choices of $r$ and $b$.
It follows that the subgraph induced by the red and blue vertices is complete bipartite.
\end{Proof}

This theorem has an application to the family of maximal outerplanar graphs.
Recall that a maximal outerplanar graph (MOP) is obtained from a cycle by adding noncrossing 
chords to make it a triangulation.
The above theorem, together with Proposition~\ref{p:generalP1} and the fact that a MOP has chromatic number $3$, shows that a
MOP with large enough diameter has \Pone-compelling chromatic number equal to~$5$.

\subsection{Graphs with large compelling number}

It is not hard to see that a graph has \Pone-compelling chromatic number $n$ if and only if it is complete. 
We conclude this section with a more general statement:

\begin{lemma}
Let $G$ be a graph of order $n$.
If $\chi(G) \ge n/2+1$, then the \Pone-compelling chromatic number of $G$ is $\chi(G)$.
\end{lemma}
\begin{Proof}
Consider a graph $G$ with $\chi(G) \ge n/2+1$.
We claim that every proper coloring is compelling. For, by averaging there must be at least two colors that 
are used only once, say on vertices $u$ and $v$. The minimality of the number of colors means that vertices $u$ and $v$ 
must be adjacent, and so every RC contains an edge.
\end{Proof}

The $n/2+1$ bound in the above lemma is improvable slightly to $(n+1)/2$. (We omit the details.) But that value is best possible because of the following example.
Let $S_m$ be the split graph obtained from an independent set $J=\{u_1, \ldots, u_m\}$ and a clique
$K=\{v_1, \ldots, v_m\}$ and joining every $u_i$ to every $v_j$ except when $i=j$. The result $S_m$
has order $2m$, chromatic number $m$, and for any coloring with $m$ colors, the set $J$ is an RC that is independent.

\section{The Connectivity-Compelling Chromatic Number}

Define the following properties of rainbow committees:

\begin{quote} 
Property~\Ptwo:
The subgraph induced by each RC is connected;\\[2mm]
Property \CDom: Every RC forms a connected dominating set of the graph.
\end{quote}

That is, Property~\CDom{} means that the subgraph induced by each RC is connected, and
every vertex not in the RC is adjacent to some vertex in the RC.

\begin{lemma}
Let $G$ be a graph with at least one edge. Then a proper coloring is  \Ptwo-compelling 
if and only if it is \CDom-compelling.
\end{lemma}
\begin{Proof}
Obviously Property~\CDom{} is at least as strong as Property \Ptwo. So we need only show the latter implies the former.
Suppose $G$ is a graph satisfying Property~\Ptwo{} but not Property~\CDom. Consider
a proper coloring of $G$ and some RC $J$ 
that fails to be a connected dominating set of $G$. By supposition, the set~$J$ is 
connected, so it must not be dominating. Thus there is some vertex $v$ that is not in $J$ and neither
is any of its neighbors. But then one can form a new RC $J'$ by taking $v$
instead of its color representative in~$J$. The graph induced by $J'$
is disconnected provided it has at least two vertices. That is,
we obtain a contradiction unless all vertices of $G$ have the same color.
\end{Proof}

It follows that if a graph $G$ is connected, for compelling purposes the properties \Ptwo{} and \CDom{} are equivalent.

By the above discussion and Proposition~\ref{p:general} we get the following bound.
Let $\conndom(G)$ denote the connected domination number of a graph $G$.
(The upper bound is equivalent to selecting a minimum connected dominating set and giving its members unique colors and giving the remainder of the graph any proper coloring.)

\begin{prop}  \label{p:generalP2}
For a nontrivial connected graph $G$, the \Ptwo-compelling chromatic number is at least $\max\{ \chi(G), \conndom(G)\}$ and at most $\chi (G)+ \conndom(G)$.
\end{prop}

\subsection{Paths, cycles, and trees}

\begin{lemma}
For $n\ge 3$, the \Ptwo-compelling chromatic number for the path $P_n$ is $ n-1$. 
\end{lemma}
\begin{Proof}
A compelling coloring is obtained by 
coloring both ends of the path the same color and giving every other vertex a unique color:
every RC induces a path.

Consider a compelling coloring.
Suppose there is some vertex $v$ that is not an end-vertex but has the same color as some other vertex $w$.
Let $T'$ be the component of the graph created by removal of $v$ that does not contain $w$.
Then consider an RC that contains $w$ and a vertex of $T'$: that RC cannot induce a connected subgraph.
That is, every interior vertex must have a unique color; thus there must be at least $n-1$ colors.
\end{Proof}

It is easy to see that the $3$-cycle has compelling number $3$ while the $4$-cycle has compelling number $2$.

\begin{lemma}
For $n\ge 5$, the \Ptwo-compelling chromatic number for the cycle $C_n$ is $n-1$.
\end{lemma}
\begin{Proof}
If one uses $n-1$ colors, then every RC will automatically form a path on $n-1$ vertices and thus be connected.
So we need only show that $n-2$ colors allows an RC that is disconnected.

Consider some compelling coloring.
Suppose some color, say red, is used three times, say at vertices $r_1$, $r_2$, and $r_3$. 
Then if one starts an RC with $r_1$ and a vertex that lies between $r_2$ and $r_3$, the resultant RC cannot 
be connected. Thus it follows that every color is used at most twice.

Suppose color red is used on vertices $r_1$ and $r_2$ and color blue on vertices
$b_1$ and~$b_2$. If these vertices do not alternate on the cycle, then we may assume
that they appear in order $r_1 r_2 b_1 b_2 r_1$ on the cycle. Then any RC
containing $r_1$ and $b_1$ will be disconnected, as it would be impossible to get from $r_1$
to $b_1$. So we may assume these vertices alternate on the cycle, say $r_1 b_1 r_2 b_2 r_1$.
Since $n\ge 5$, there is at least one other vertex, say $v$ (necessarily of another color). Say $v$ lies
between $r_1$ and $b_1$. Then any RC containing the vertices $v$, $r_2$, and $b_2$ is disconnected.

So we have shown that at most one color is repeated, and it only appears twice. Thus $n-1$ colors are needed.
\end{Proof}

The above result for the path generalizes to trees. Let $\int(T)$ denote the number of interior (non-leaf) vertices
of the tree $T$.

\begin{lemma}
For any tree $T$ of order at least $3$, the \Ptwo-compelling chromatic number equals $1+\int(T)$.
\end{lemma}
\begin{Proof}
For a compelling coloring of $T$, give all leaves the same color and give each interior vertex a unique color.

Consider a compelling coloring.
Suppose some interior vertex $v$ has the same color as some other vertex $w$.
Let $T'$ be any component of $T-v$ not containing~$w$. Then if one starts an RC with 
$w$ and a vertex in $T'$, the resultant RC cannot be connected. It follows that every interior vertex must have its own unique color.
\end{Proof}

\subsection{Extremal values}

The parameter equals $2$ if and only if the graph is complete bipartite. For,
to equal $2$, one needs a graph that is bipartite; if one takes an RC the two vertices must be adjacent: so the graph must be
complete bipartite.

The parameter equals $n$ if and only if complete. For, if not complete, it is possible to find
two vertices $u$ and $v$ that are not adjacent and neither is a cut-vertex. Give them the same color
and color each other vertex with a unique colors.

One obvious question is whether large chromatic number implies that the \Ptwo-compelling chromatic number is the chromatic number.
The answer is no, as even a graph with chromatic number $n-2$ can have larger \Ptwo-compelling chromatic number. For example,
take the complete graph $K_{n-1}$, remove one edge $e$, and then add a leaf adjacent to one end of $e$. The result
is the graph with chromatic number $n-2$ and \Ptwo-compelling number $n-1$ shown below.

\begin{center}
\includegraphics{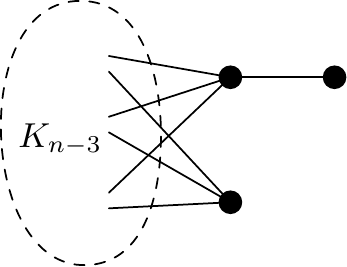}
\end{center}

\subsection{MOPs}

We conclude this section by showing that the \Ptwo-compelling chromatic number of a maximal outerplanar graph is determined by
its connected domination number. By Proposition~\ref{p:generalP2} we know that the parameter is at 
least $\conndom(G)$ and at most $\conndom(G)+3$. But:

\begin{theorem}
\label{t:mopP2}
Let $G$ be a MOP on at least three vertices. Then
the \Ptwo-compelling chromatic number of $G$ is $\conndom(G)+2$.
\end{theorem}

We first prove the upper bound.

\begin{claim}  \label{c:one}
If one removes a connected dominating set $D$ from an outerplanar graph $G$,
then what remains has no cycle. 
\end{claim}
\begin{Proof}
It is known that outerplanar graphs are closed under contractions and that $K_4$ is not outerplanar.
Suppose there were a cycle~$C$ in $G-D$. Then when one contracts $D$ to a single vertex $d$, vertex $d$ would be adjacent to all of $C$ and so
one would have a subdivision of $K_4$.
\end{Proof}

It follows that one can produce a compelling coloring of $G$ by taking a minimum connected dominating set $D$, giving each of its vertices a unique color,
and then using two colors on the remaining forest. This coloring has $\conndom(G)+2$ colors.

We now prove the lower bound. 
The theorem is true for a MOP with three vertices; so assume $G$ has at least four vertices. Thus there is at least one chord.
For a coloring, define $\un$ as the set of vertices that receive unique colors.
Define a \emph{chord-cover} of $G$ as a set of vertices that intersects every chord of $G$.
For example, in the following MOP the white vertices form a (minimal) chord-cover.

\begin{center}
\includegraphics{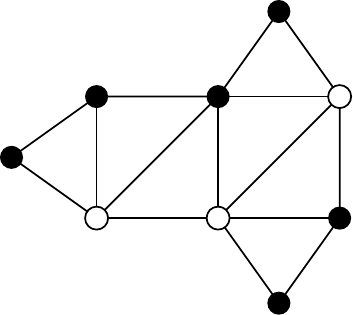}
\end{center}

\begin{claim} \label{c:two}
If a coloring of a MOP is \Ptwo-compelling, then $\un$ forms a chord-cover.
\end{claim}
\begin{Proof}
Suppose there is some chord $xy$ neither end of which is in $\un$. Then build an RC starting with another vertex of the color of $x$ 
and another vertex of the color of $y$. Note that the removal of the two vertices $x$ and $y$ splits the MOP $G$ into two components. Say the RC
already has a vertex from the first component.
Then $xy$ lies in a triangle with some vertex $z$ in the second component, and vertex $z$ has a different color and so can be added to the nascent RC. 
The resulting RC will be disconnected, a contradiction. Thus it is necessary that $\un$ forms a chord-cover. 
\end{Proof}

\begin{claim} \label{c:three}
In any MOP with at least $4$ vertices, a chord-cover $C$ is a connected dominating set.
\end{claim}
\begin{Proof}
For domination, note that every vertex $v$ is in some triangle and
every triangle uses at least one chord; thus $v$ is in $C$ or has a neighbor in $C$.

For connectivity, given any two vertices $c$ and $d$ in $C$ there is sequence of triangles $T_1, \ldots, T_k$ 
with $c$ in $T_1$ and $d$ in $T_k$
such that consecutive triangles overlap
in a chord. Assume one has a shortest such sequence. We prove by induction on~$k$ that $c$ and $d$ are
connected using only vertices of $C$. If $k=1$ then $c$ and $d$ are adjacent. So assume $k\ge 2$.
Then by minimality, vertex $c$ is not in~$T_2$. 
The overlap between $T_1$ and $T_2$ is a chord, and thus one of those vertices, say~$f$, must be in~$C$.
By the inductive hypothesis vertices $f$ and $d$ are connected inside~$C$, while vertices $f$ and $c$ are adjacent. Thus
$c$ and $d$ are connected using only vertices of~$C$.
\end{Proof}

\begin{claim} \label{c:four}
In any MOP with at least $3$ vertices, 
the complement of a minimal connected dominating set~$D$ must contain an edge.
\end{claim}
\begin{Proof}
Suppose the complement of $D$ is an independent set. Then every vertex outside $D$ has at least two neighbors in $D$.
Let $d$ be a non-cut-vertex of the subgraph induced by $D$. (Every graph has such a vertex.) Then the set $D-\{d\}$ is still connected
and still dominating, a contradiction.
\end{Proof}

So we return to proving the lower bound for Theorem~\ref{t:mopP2}
 and consider a compelling coloring.
By Claims~\ref{c:two} and~\ref{c:three} we know that $\un$ must be a connected dominating set. If $\un$
has size exactly $\conndom(G)$, then by Claim~\ref{c:four} (at least) two extra colors are needed for the remaining vertices, for a total of $\conndom+2$ colors.
Further, if $\un$ has size exactly $\conndom(G)+1$, then it cannot be the whole graph, and so one needs at least one
color for what remains. That is, a compelling coloring uses at least $\conndom(G)+2$ colors.
And thus the  theorem is true.

\section{Complexity}

We conclude with some thoughts on the complexity of the parameters.
As observed by~\cite{GHR06,HHLMW09}, for fixed $k\ge 4$ the problem of determining whether the dominator chromatic number is 
at most $k$ is NP-hard. The idea is to add a vertex adjacent 
to all vertices: the resultant graph $G'$ has a dominator coloring with $4$ colors if and only
if the original graph $G$ is $3$-colorable (which is an NP-hard problem). 
The result was extended to total dominator chromatic number in~\cite{HHHMR,Kazemi15}.
Presumably the generalization is NP-hard for $k\ge 4$ ``usually''. For example, the join-of-a-dominator-vertex 
construction shows that determining whether the \Pone-compelling or \Ptwo-compelling chromatic number is at most $k$ is NP-hard for $k\ge 4$.

In contrast, testing whether the dominator chromatic number is $3$ 
was shown in \cite{HHHMR} to have a polynomial-time algorithm.
We show here that 
testing whether a graph has total dominator chromatic number~$3$ can also be performed in polynomial time.
Since a subset of size $3$ is a total dominating set if and only if it is a connected dominating set, 
a consequence is that testing whether the \Ptwo-compelling chromatic number is $3$ is polynomial-time computable too.

\begin{theorem}
Testing whether a graph has total dominator chromatic number $3$ can be performed in 
polynomial time.
\end{theorem}

\begin{Proof}
The requirement of a total dominator coloring is that for each vertex $v$ there is a color $c_v$ different
from the color of $v$ such that $v$ is adjacent to all vertices of color $c_v$.
One has to determine whether such a coloring exists. There are two possibilities for such a colorig.
\smallskip

\textbf{Case 1:} \textsl{For each color, all vertices $v$ of that color have the same~$c_v$.}
Say every red vertex $v$ has $c_v=$ blue. This means that every red vertex is joined to every blue vertex.
This immediately implies that every blue vertex $v$ has $c_v=$ red. 
And then say every green vertex $v$ has $c_v=$ red. It follows that every red vertex has the same neighborhood.
Conversely, a coloring where every non-red-vertex is adjacent to every red vertex is a total dominator coloring and has the attributes for this case.
Thus one can test for such a situation by guessing  a red vertex~$v$ (that is, trying all possibilities), checking that it 
and each of its nonneighbors are clones, then giving it
and each of its nonneighbors red, and finally testing whether what remains is bipartite.\smallskip

\textbf{Case 2:} 
\textsl{There are two vertices $u$ and $v$ of the same color, say red, with $c_u \neq c_v$.} 
Say $u$ is adjacent to all blue vertices (but not all greens) and $v$ to all green vertices (but not all blues). 
In particular that means that all the blue and green vertices are contained within $N(u) \cup N(v)$ and 
indeed are exactly $N(u) \cup N(v)$, since none of those vertices can be red. And thus the red vertices are precisely the vertices
outside $N(u) \cup N(v)$.\smallskip

\textbf{Case 2.1:}
\textsl{There are two vertices $x$ and $y$ of another color, say blue, with $c_x \neq c_y$.} 
By the same reasoning this
determines the blue vertices. And thus the green vertices are determined. So one can test for such a situation
by guessing the vertices $u,v,x,y$ (that is, trying all possibilities) and checking the resultant coloring for the desired property.\smallskip

\textbf{Case 2.2:}
\textsl{Every blue vertex $x$ has the same $c_x$ and every green vertex $y$ has the same $c_y$.}
Since there is some blue vertex not adjacent to $v$ and some green vertex not adjacent to $u$, 
it must be that every blue vertex is adjacent to every green vertex (and necessarily vice versa). One can
test for this situation by guessing the vertices $u$ and $v$ and checking whether $N(u) \cup N(v)$ induces a complete bipartite subgraph, thereby
yielding the coloring, and then checking whether each red vertex has the desired property.
 \end{Proof}

\section{Final Thoughts}

We conclude with some thoughts about further work. One natural idea is to continue studying the general schema. As regards the two specific parameters introduced here, their values for cubic graphs 
and for disjoint unions seem interesting. Further,
the cases that the property of the RC is that it 
contains at least two edges, or that it contains a cycle, seem worthy of study.

\end{document}